\documentclass[11pt,oneside]{article}
\usepackage{amssymb,amsmath}
\usepackage{url}
\renewcommand{\le}{\leqslant}
\renewcommand{\ge}{\geqslant}

\newcommand{\rad}{\mathrm{rad}}

\newcommand{\RR}{\mathbb{R}}
\newcommand{\ZZ}{\mathbb{Z}}
\newcommand{\QQ}{\mathbb{Q}}

\newcommand{\NN}{\mathbb{N}}
\newcommand{\FF}{\mathbb{F}}

\newcommand{\KK}{\mathbb{K}}

\renewcommand{\AA}{\mathbb{A}}

\newcommand{\va}{\mathbf{a}}

\newcommand{\sep}{\mathrm{Sep}}

\newtheorem{theorem}{Theorem}

\newtheorem{problem}{Problem}

\newtheorem{conjecture}[problem]{Conjecture}

\newtheorem{corollary}{Corollary}
\newtheorem{theoremrl}{Theorem RL}

\newcommand{\bigk}{\mathop{\mathbf{K}}}

\begin{document}
\title{Distance between cubics and rationals}

\author{Dmitry Badziahin}

%\thanks{The University of Sydney, dzmitry.badiahin@sydney.edu.au}

\maketitle

\begin{abstract}
We investigate the following problem: what is the smallest possible
distance between a cubic irrational $\xi$ and a rational number
$p/q$ in terms of the height $H(\xi)$ and $q$? More precisely, we
consider the set $D_{3,1}$ consisting of all pairs $(u,v)$ of
positive real numbers such that $|\xi - p/q| > cH^{-u}(\xi)q^{-v}$
for all cubic irrationals $\xi$ and rationals $p/q$. First, we
transform this problem into one about the root separation of cubic
polynomials. Second, under the assumption of the famous
abc-conjecture, we give an almost complete description of $D_{3,1}$.
Namely, the points $(u,v)$ with $2\le v\le 3$ that lie in the
interior of $D_{3,1}$ are characterised by the inequality $u>
10-3v$. Assuming only the weaker Hall conjecture, we also obtain
nontrivial results about the shape of $D_{3,1}$, although these are
not as strong as those derived from the abc-conjecture. Finally, we
discuss an analogue of the set $D_{3,1}$ in function fields where we
are able to give a complete description unconditionally.
\end{abstract}

{\footnotesize{{\em Keywords}: Roth theorem, approximation to
algebraic numbers, continued fractions, root separation.

Math Subject Classification 2020: 11J68, 11D59}}

\section{Introduction}

For a positive integer $d$, let $\AA_d$ denote the set of real
algebraic numbers of degree $d$. We equip each element $\xi\in
\AA_d$ with its naive height $H(\xi)$, defined as the maximum
absolute value of the coefficients of its minimal polynomial
$P_\xi$. In the classical problem from Diophantine approximation,
one fixes an algebraic number $\xi$ and studies the rate of its
approximation by rational numbers $p/q$ or more generally, by
numbers $\alpha$ from a fixed number field $\KK$. Thanks to the
landmark works of Roth and LeVeque~\cite{roth_1955, leveque_1956},
this rate is well understood.

\begin{theoremrl}
Let $K$ be an algebraic number field. For any real algebraic
$\xi\not\in K$ and $\lambda>2$ the equation
$$
|\xi - \alpha|< H(\alpha)^{-\lambda}
$$
has only finitely many solutions $\alpha\in K$.
\end{theoremrl}

In view of the classical Dirichlet theorem, the lower bound on the
degree $\lambda$ here is best possible. It is widely believed that
algebraic numbers of degree greater than 2 are not badly
approximable, and therefore Theorem~RL is thought not to hold for
$\lambda=2$. However, this remains an open problem.

In this paper, we introduce a different though related problem.
Given two distinct positive integers $d$ and $r$, we investigate
functions $\Psi_{d,r}:\NN^2\to \RR^+$ for which the inequality
\begin{equation}\label{main_eq}
|\xi-\alpha| \ge \Psi_{d,r}(H(\xi), H(\alpha))
\end{equation}
holds for all pairs $(\xi,\alpha)\in \AA_d\times \AA_r$.
Alternatively, one may pose a slightly weaker version of the
question by requiring~\eqref{main_eq} to hold for all but finitely
many such pairs. If for some absolute constant $c_{d,r}$ the
function $\Psi_{d,r}(x,y) = c_{d,r}x^{-u}y^{-v}$
satisfies~\eqref{main_eq} we say that $\AA_d$ and $\AA_r$ are
$(u,v)$-distanced. One can easily check that this property does not
depend on whether we require~\eqref{main_eq} to hold for all pairs
or for all but finitely many pairs $(\xi,\alpha)\in \AA_d\times
\AA_r$. Moreover, it is immediate that if $\AA_d$ and $\AA_r$ are
$(u_0,v_0)$-distanced then they are also $(u,v)$-distanced for all
$u\ge u_0$ and $v\ge v_0$. Formally, the main problem of this paper
is:

\begin{problem}
Given two distinct positive integers $d,r$, determine the set
$$
D_{d,r}:=\{ (u,v)\in (\RR^+)^2: \AA_d \mbox{ and } \AA_r \mbox{ are
} (u,v)\mbox{-distanced}\}.
$$
\end{problem}

The classical Liouville inequality provides one element of
$D_{d,r}$: namely, $$\Psi_{d,r}(x,y) = (d+1)^{-r}(r+1)^{-d}
x^{-r}y^{-d}$$ satisfies condition~\eqref{main_eq}. In other words,
$\AA_d$ and $\AA_r$ are always $(r,d)$-distanced. If $r=1$, i.e.
$\AA_r=\QQ$, the Dirichlet theorem implies that every point
$(u,v)\in D_{d,1}$ must satisfy $v\ge 2$. On the other hand, the
folklore conjecture claims that for $d>2$ all numbers in $A_d$ are
not badly approximable. If true, this would imply that no pairs of
the form $(u,2)$ belong to $D_{d,1}$ as soon as $d>2$. By fixing
$p/q\in\QQ$ and carefully choosing an infinite sequence of
irreducible polynomials $P\in \ZZ[x]$ such that $q^3P(p/q)=1$, one
can also deduce that every element $(u,v)\in D_{d,1}$ satisfies
$u\ge 1$ (though we do not rigorously prove this claim in the
present paper). Taken together, these facts show that the only
non-trivial part of $D_{d,1}$ lies in the range $2\le v<d$, which is
non-empty as soon as $d>2$.

Let $S_d:=\{(x,y)\in \RR^2: 2\le y<d\}$ and $L_u:= \{(x,y)\in \RR^2:
y=u\}$. So far, the shape of $D_{d,1}\cap S_d$ for $d\ge 3$ remained
completely mysterious. We do not even have a reasonable conjecture
about it. Establishing that $D_{d,1}\cap L_v \neq \emptyset$ for
some $v<d$ would represent a major step toward an effective
improvement of Liouville's theorem which still remains out of reach.

Here we focus on the smallest non-trivial set $D_{3,1}$. Our first
step is to convert the problem of determining the shape of $D_{3,1}$
into a problem concerning the root separation in cubic polynomials.
Recall that the root separation $\sep(Q)$ of a polynomial
$Q\in\RR[x]$ of degree at least 2 is defined as the distance between
its two closest roots.

\begin{theorem}\label{th1}
Let $Q\in\ZZ[x]$ be an irreducible cubic polynomial. Denote by $B$
the absolute value of its leading coefficient and set $A:= H(Q)/B$,
where $H(Q)$ is the naive height of $Q$. Suppose that
\begin{equation}\label{th1_eq}
\sep(Q) \ge \frac{c}{B^{2+s}A^{2-t}}
\end{equation}
holds for all such $Q$ and for some positive constants $c,s,t$. Then
$$
(u,v) = \left(\frac{2+s}{s+t}, 2+\frac{s}{s+t}\right) \in D_{3,1}.
$$
\end{theorem}

In this paper we adopt Vinogradov's notation. For positive
quantities $A$ and $B$ we write $A\gg B$ if $A\ge cB$ for some
absolute constant $c>0$. If the constant $c$ is allowed to depend on
the parameter $t$, we write $a\gg_t B$. Similarly, we define $A\ll
B$ and $A\ll_t B$. Finally, we write $A\asymp B$ (respectively,
$A\asymp_t B$) if both $A\ge B$ and $A\ll B$ hold (respectively,
$A\gg_t B$ and $A\ll_t B$).

It is known, due to Mahler~\cite{mahler_1964}, that the root
separation satisfies $\sep(Q)\gg H^{-2}(Q) = (AB)^{-2}$. However,
very little is known beyond this, in particular regarding lower
bounds of the form~\eqref{th1_eq}. While the existence of positive
parameters $s$ and $t$ satisfying~\eqref{th1_eq} seems plausible,
this condition has not been extensively studied. The only related
work we can find is by Bugeaud and Mignotte~\cite{bug_mig_2010}, who
consider the extremal case $B=1$ of this problem.

While verifying~\eqref{th1_eq} is likely to be difficult,
conditional results in this direction can be obtained. It is well
known that an effective version of Theorem~RL follows from an
effective version of the classical abc-conjecture, see for
example~\cite[Chapter 12.2]{bom_gub_2006}. Moreover, the machinery
developed there allows the construction of functions $\Psi_{d,1}$ of
the form $\Psi_{d,1}(x,y) = c\psi(x)y^{-2-\epsilon}$ that
satisfy~\eqref{main_eq}. However, the component $\psi(x)$ decays
very rapidly and can not be presented in the form $\psi(x) =
x^{-u}$. In Section~\ref{sec3}, we provide conditional non-trivial
lower bounds on $\sep(Q)$ of the form~\eqref{th1_eq}. The weaker one
assumes the Hall conjecture, while the stronger one assumes the more
general abc-conjecture. For convenience of the reader we formulate
both of them below.

\begin{conjecture}[Hall]
For any $\epsilon>0$ the inequality $$ 0<|x^3 -
y^2|<|x|^{1/2-\epsilon}$$ has only finitely many solutions.
\end{conjecture}

\begin{conjecture}[abc]
Let $a,b,c$ be coprime integer numbers such that $a+b=c$. Then for
any $\delta>0$ one has
$$
\max\{|a|, |b|, |c|\}\ll_\delta \rad(abc)^{1+\delta},
$$
where $\rad(n)$ is the radical of $n$, i.e. it is the product of all
distinct prime numbers that divide~$n$.
\end{conjecture}

Under the condition of the abc-conjecture, one can provide an almost
complete description of the set $D_{3,1}$, which is done in the
following

\begin{theorem}\label{th5}
Suppose that the abc-conjecture is satisfied. Then
$$
\{(u,v)\in\RR^2: 2<v\le 3, u>10-3v\} \subset D_{3,1}.
$$
On the other hand, the following inclusion holds unconditionally:
$$
D_{3,1}\cap S_3\subset\{(u,v)\in \RR^2: 2<v\le 3, u\ge 10-3v\}.
$$
\end{theorem}

This result provides a precise description of the interior of the
set $D_{3,1}\cap S_3$. What remains unknown is which points on the
boundary of $D_{3,1}$ belong to this set. The first part of
Theorem~\ref{th5} follows from the following result, proved in
Section~\ref{sec3}:

\begin{theorem}\label{th2}
Suppose that the Hall conjecture is satisfied for some $\epsilon>0$.
Then for all irreducible cubic polynomials $Q$ one has
\begin{equation}\label{th2_eq}
\sep(Q)\gg_\epsilon \frac{1}{B^{2+\frac12 - \epsilon}A^{2-\frac12
-\epsilon}}.
\end{equation}
Here $B$ is the absolute value of the leading coefficient of $Q$ and
$A = H(Q)/B$. This in turn implies that for all $\frac12+\epsilon\le
r\le 1$,
\begin{equation}\label{th2_eq2}
(u,v) = \left(r + \frac{2(1-r)}{1/2-\epsilon}, 2+r\right) \in
D_{3,1}.
\end{equation}
If one assumes the more general abc-conjecture for $\delta =
\frac{\epsilon/3}{1-\epsilon/3}$ then the polynomials $Q$ satisfy
$$
\sep(Q)\gg_\epsilon \frac{1}{B^{2-\epsilon} A^{2-\frac12
-\epsilon}},
$$
and the inclusion~\eqref{th2_eq2} holds for all $2\epsilon\le r\le
1$.

Moreover, one can explicitly compute the constant $c$ such that the
inequality
\begin{equation}\label{th2_eq3}
\left|\xi-\frac{p}{q}\right| > \frac{c}{H^u(\xi) q^v}
\end{equation}
is satisfied for all but finitely many pairs $(\xi,p/q)\in
\AA_3\times\AA_1$. If an effective version of the Hall conjecture
(respectively, the abc-conjecture) is true then the last inequality
can be assumed for all pairs $(\xi,p/q)$.
\end{theorem}

As a corollary of Theorem~\ref{th2}, one can establish the following
result about Thue equations. It is done in Section~\ref{sec6}.
\begin{corollary}\label{corl}
For $\va = (a_0,a_1,a_2,a_3)\in\ZZ^4$ define the cubic binary form
$$
F_\va(x,y):= a_3x^3 + a_2x^2y + a_1xy^2 + a_0y^3.
$$
Suppose that the Hall conjecture is true for a given $\epsilon>0$.
Then the inequality
$$
0<|F_\va(p,q)|<\frac{q^{1/2-\epsilon}}{||\va||_\infty^{4+2\epsilon}}
$$
has only finitely many solutions in $(\va,p,q)\in
\ZZ^4\times\ZZ\times\ZZ$. Moreover, if an effective version of the
Hall conjecture is assumed then there exists an effectively
computable constant $c>0$ such that
$$
F_\va(p,q)=0\quad\mbox{or}\quad
|F_\va(p,q)|>\frac{cq^{1/2-\epsilon}}{||\va||_\infty^{4+2\epsilon}}
$$
for all integer $(\va,p,q)$.

If the abc-conjecture is true for $\delta =
\frac{\epsilon/3}{1-\epsilon/3}$ where $\epsilon>0$ is given, then
the inequality
$$
0<|F_\va(p,q)|< \frac{q^{1-2\epsilon}}{||\va||_\infty^{5+4\epsilon}}
$$
has only finitely many solutions in $(\va,p,q)\in
\ZZ^4\times\ZZ\times\ZZ$.
\end{corollary}

In Section~\ref{sec4}, we present an infinite family of pairs
$(\xi,p/q)\in\AA_3\times \QQ$ that sit very close to each other.
This will allow us to verify the second part of Theorem~\ref{th5}.
We formulate it separately below.
\begin{theorem}\label{th3}
The point $(u,v)\in D_{3,1}$ with $2\le v\le 3$ must satisfy $u\ge
10-3v$.
\end{theorem}

The pairs $(\xi,p/q)\in\AA_3\times \QQ$ from Section~\ref{sec4} look
quite mysterious. We leave the construction of such examples as an
interesting and challenging exercise for the curious reader. We
verify that the continued fractions of these cubic $\xi$ start with
a very unusual pattern. While we can confirm this pattern
numerically, the underlying reason for its occurrence remains
unknown. We briefly discuss this phenomenon in Section~\ref{sec4}.

We conclude this paper with a discussion of the analogue of the set
$D_{3,1}$ over function fields. Given a base field $\FF$, we define
the set of formal Laurent series with coefficients in $\FF$ by
$\KK:=\FF((t^{-1}))$. Like in the real case, one can measure the
distance between algebraic Laurent series of degrees $d$ and $r$ and
define an analogue $D_{d,r}(\KK)$ of the set $D_{d,r}$. A precise
definition and further discussion of these sets can be found in
Section~\ref{sec5}. Unlike the case of real numbers, we can describe
the shape of the set $D_{3,1}(\KK)$ unconditionally, thanks to the
works of Osgood~\cite{osgood_1973, osgood_1975}. Notably, the shape
of $D_{3,1}(\KK)$ depends crucially on the characteristics of the
base field $\FF$.

\begin{theorem}\label{th4}
If $\FF$ has positive characteristic $p$ then for all $d\ge 2$,
$$
D_{d,1}(\KK) = \{(u,v)\in\RR: u\ge 1,v\ge d\}.
$$
Otherwise, if characteristic of $\FF$ is zero then for all $2\le
v\le 3$ the set $D_{3,1}(\KK)\cap L_v$ is nonempty. Moreover,
$(4,2)\in D_{3,1}(\KK)$.
\end{theorem}
Since our main focus is on the set $D_{3,1}$ in the real setting, we
do not investigate the precise shape of $D_{d,1}(\KK)$, $d\ge 4$ for
the base fields $\FF$ of characteristic zero. But that would be an
interesting direction for further research.

%With help of computational experiments and various heuristics we
%will try to understand the shape of $D_{3,1}\cap S_3$.

\section{Relation with the root separation problem}\label{sec2}

Consider a cubic irrational $\xi\in\RR$. For convenience, assume
that $0\le \xi\le 1$. Let its minimal polynomial be
$$
P_\xi:=a_3x^3+a_2x^2+a_1x+a_0.
$$
Denote the other two zeroes of $P_\xi$ by $\xi_2, \xi_3$. Suppose
there is a rational number $p/q$ such that
\begin{equation}\label{eq9}
\left|\xi - \frac{p}{q}\right| \le \frac{c_1}{H(\xi)^u q^v}
\end{equation}
for some $u\ge 2, v>2$ and small absolute constant $c_1>0$. Denote
the left hand side of this expression by $\frac{1}{q^2A}$. Recall
that by the separation of a polynomial $P$ without double zeroes we
call the minimal distance between its roots. The classical result of
Mahler states that
$$
\sep(P_\xi)\gg \frac{1}{H(\xi)^2}.
$$
Therefore, as soon as we have $u\ge 2$, $v>0$ and small enough
constant $c_1$, the following inequalities are satisfied:
\begin{equation}\label{eq1}
(1-\epsilon)|\xi_i-\xi|\le \left|\xi_i-\frac{p}{q}\right|\le
(1+\epsilon)|\xi_i-\xi|,\quad i=2,3,
\end{equation}
where $\epsilon>0$ can be chosen arbitrarily small.

Since $p/q$ is too close to $\xi$, by the Legendre theorem it has to
be the convergent of $\xi$: $p/q = p_n/q_n$. Let $p_{n+1}/q_{n+1}$
be the next convergent of $\xi$ and consider the polynomial $Q$ that
has roots
$$
\eta_1:=\frac{q_{n+1}\xi-p_{n+1}}{q_n\xi - p_n},\quad \eta_2:=\frac{q_{n+1}\xi_2-p_{n+1}}{q_n\xi_2-p_n},\quad \eta_3:=\frac{q_{n+1}\xi_3-p_{n+1}}{q_n\xi_3-p_n}.
$$
From basic properties of convergents, we have that
$q_{n+1}p_n-p_{n+1}q_n = \pm 1$, therefore the discriminants of
$P_\xi$ and $Q$ coincide. One can check that $Q$ can be written in
the form
$$
a_3(p_nx-p_{n+1})^3 + a_2(p_nx-p_{n+1})^2(q_nx-q_{n+1}) + a_1(p_nx-p_{n+1})(q_nx-q_{n+1})^2 + a_0 (q_nx-q_{n+1})^3.
$$
Indeed, from the definitions of $\eta_i$, $i\in\{1,2,3\}$, one
derives that $$\xi_i = \frac{p_n\eta_i-p_{n+1}}{q_n\eta_i -
q_{n+1}}.$$ Therefore the numbers $\eta_1,\eta_2$ and $\eta_3$ are
roots of the polynomial
$(q_nx-q_{n+1})^3P_\xi\left(\frac{p_nx-p_{n+1}}{q_nx-q_{n+1}}\right)$.
In particular, the leading coefficient of $Q$ equals
\begin{equation}\label{eq11}
|b_3| = q^3 |P_\xi(p/q)| = \frac{q}{A}
(1+\delta)|a_3(\xi-\xi_2)(\xi-\xi_3)|
\end{equation}
where $\delta$ is some small parameter, $|\delta|\le
2\epsilon+\epsilon^2$. Notice that $|P'_\xi(\xi)| = |3a_3\xi^2 +
2a_2\xi + a_1|\le 6H(\xi)$. On the other hand, it equals
$|a_3(\xi-\xi_2)(\xi-\xi_3)|$. Therefore we get an upper bound
\begin{equation}\label{eq2}
|b_3|\le \frac{7 q H(\xi)}{A}.
\end{equation}

Now we provide the bounds for the roots of $Q$.
$$
\left|\frac{q_{n+1}\xi-p_{n+1}}{q_n\xi-p_n}\right| < 1,
$$
by the choice of convergents $p_n/q_n, p_{n+1}/q_{n+1}$ of $\xi$. On
the other hand, notice that
$$
\frac{1}{q_{n+1}+q_n}< |q_n\xi-p_n| < \frac{1}{q_{n+1}}.
$$
Therefore $q_{n+1}$ is between $(A-1)q_n$ and $Aq_n$. By choosing
the constant $c_1$ in the assumption~\eqref{eq9} small enough we
guarantee that $A$ is large enough, and thus $q_{n+1} = Aq
(1+\delta_1)$ where $|\delta_1|\le \epsilon$. We are now ready to
estimate the other roots:
$$
\left|\frac{q_{n+1}\xi_i-p_{n+1}}{q_n\xi_i-p_n}\right| =(1+\delta_2)
\frac{q_{n+1}}{q_n} \cdot \frac{|\xi-\xi_i|}{|\xi-\xi_i|} =
(1+\delta_3)A,
$$
where $i=2,3$ and $|\delta_3|\le 2\epsilon+\epsilon^2$.

Next, we compute the distance between the roots
\begin{equation}\label{sepeq}
\left|\frac{q_{n+1}\xi_2-p_{n+1}}{q_n\xi_2-p_n} -
\frac{q_{n+1}\xi_3-p_{n+1}}{q_n\xi_3-p_n}\right|
\stackrel{\eqref{eq11}}=
\frac{|(1+\delta)a_3^2(\xi-\xi_2)(\xi-\xi_3)(\xi_2-\xi_3)|}{b_3^2A^2}\asymp
\frac{\sqrt{\Delta(P_\xi)}}{b_3^2A^2},
\end{equation}
where $\Delta(P_\xi)$ is the discriminant of $P_\xi$ which satisfies
$|\Delta(P_\xi)|\le 54H^4(\xi)\ll H^4(\xi)$.

The upshot is that if $\xi$ and $p/q$ are too close to each other
then there exists another polynomial $Q$ whose two of three roots
are very close to each other. To simplify it a little bit more, we
consider a polynomial $R$ whose roots are $\sigma_i = \eta_i-k$
where $k\in\ZZ$ is such that the real part of $\sigma_2$ is between
-1/2 and 1/2. Then the real part of $\sigma_3$ is between
$-1/2-\epsilon$ and $1/2+\epsilon$, and $|\sigma_1| = (1+\delta_4)A
\asymp A$ where $|\delta_4|\le 3\epsilon+\epsilon^2$. Such a shift
does not change the leading coefficient of the polynomial therefore
the one for $R$ is $b_3$. On the other hand,
$$
|b_2|/|b_3| = |\sigma_1+\sigma_2+\sigma_3| = (1+\delta_5)A \asymp
A,\; |\delta_5| \le 4\epsilon+\epsilon^2
$$
and by the same reason, $|b_1|\le (1+6\epsilon)A|b_3|$, $|b_0|\le
\left(\frac14+3\epsilon\right)A|b_3|$.

Recall that $A \ge q^{v-2}H^u(\xi)/c_1$ and for convenience denote
$\tau = v-2$. This equality together with~\eqref{eq2} implies that
$$
H^{u-\tau} \le 7^\tau
c_1\frac{A^{1-\tau}}{|b_3|^\tau}\quad\Longrightarrow\quad
\Delta(R)\le 54\cdot(7^\tau c_1)^{\frac{4}{u-\tau}}
A^{\frac{4(1-\tau)}{u-\tau}} |b_3|^{-\frac{4\tau}{u-\tau}}.
$$

We choose $c_1$ small enough so that the constant term on the right
hand side of the above upper bound on $\Delta(R)$ is at most $c^2$
for a given constant $c$. Then from~\eqref{sepeq}, the separation of
$R$ can be estimated as
$$
\sep(R)= \frac{(1+\delta)\sqrt{\Delta(R)}}{|b_3|^2 A^2} \le
\frac{c}{|b_3|^{2+\frac{2\tau}{u-\tau}}A^{2-\frac{2(1-\tau)}{u-\tau}}}
$$
For convenience, we denote $\rho:= \frac{2}{u-\tau}$.

The conclusion is that if the sets $\AA_3$ and $\AA_1=\QQ$ are {\bf
not} $(u,v)$-separated with $v>2$ then there exists a polynomial
$R$, whose leading coefficient $b_3$ satisfies $|b_3|=B$, the height
is $H(R) \asymp AB$ and which satisfies
$$
\sep(R) \le c B^{-2-\rho\tau}A^{-2+\rho(1-\tau)}.
$$
By adjusting $A$ and changing the absolute constant $c$, we may
assume without loss of generality that the same inequality on
$\sep(R)$ is satisfied under the condition $H(R) = AB$.

Conversely, suppose that one can show that, under the same
conditions for the polynomial $R$, one always has
$$
\sep(R) > cB^{-2-s}A^{-2+t}
$$
for some positive $c,s,t>0$. Then immediately we have that $\AA_3$
and $\AA_1$ are $(u,v)$-separated where
$$
\rho\tau = \frac{2\tau}{u-\tau} = s \quad\mbox{and}\quad
\rho(1-\tau) = \frac{2(1-\tau)}{u-\tau} = t.
$$
Solving this system of equations and noticing that $v = \tau+2$, we
derive that $D_{3,1}$ contains the pair
\begin{equation}\label{eq3}
(u,v) = \left(\frac{2+s}{s+t}, 2+\frac{s}{s+t}\right).
\end{equation}

Hence we prove Theorem~\ref{th1}. Notice that, while we do not
provide the precise formula for the relation between constants $c$
and $c_1$ here, the computations from the proof allow to construct
it easily.

%Condition of the form~\eqref{th1_eq} looks plausible, however it is
%not clear at all, what values of $s$ and $t$ should satisfy it.
%In~\cite{bug_mig_2010} the authors consider monic cubic polynomials,
%that correspond to the extremal case of $B=1$. They conclude that,
%under the conjecture of Hall, about small values of the expression
%$|x^3-y^2|$, which in turn follows from the abc conjecture, the
%value $t$ can be taken close to $1/2$.

\section{Conditional result}\label{sec3}

%It is well known that the effective version of Theorem~RL follows
%from the strong abc-conjecture, which is essentially an effective
%version of the classical abc-conjecture. See for
%example~\cite[Chapter 12.2]{bom_gub_2006}. Moreover, the machinery
%in that book allows to construct the functions $\Psi_{d,1}$ of the
%form $\Psi_{d,1}(x,y) = c\psi(x)y^{-2-\epsilon}$ that
%satisfy~\eqref{main_eq}. However, they component $\psi(x)$ decays
%very quickly and can not be presented in the form $\psi(x) =
%x^{-u}$. However, the method that links the problem about $D_{3,1}$
%and the one about the root separation in cubic polynomials, allows
%to derive a non-trivial information about this set under the
%condition of the Hall conjecture which is a particular case of the
%abc-conjecture.

Recall that in the previous section we have constructed a polynomial
$R(x) = b_3x^3+b_2x^2+b_1x+b_0\in\ZZ[x]$ such that it satisfies the
conditions
$$
|b_3| = B,\quad H(R)= AB,\quad |\eta_2-\eta_3|\ll
\frac{\sqrt{\Delta(R)}}{A^2B^2},\quad |\eta_1|\asymp A,\quad
|\eta_2|,|\eta_3|\ll 1.
$$
where $\eta_1,\eta_2,\eta_3$ are the roots of $R$ and all the
implied constants can be made explicit. The aim now is to prove
Theorem~\ref{th2}. We will do this by providing a lower bound on its
discriminant $\Delta(R)$ in terms of $A$ and $B$.

Here we follow the ideas from~\cite{bug_mig_2010}. First of all, we
simplify the discriminant by making the coefficient $b_2$ of $R$
zero. To do that, we build another polynomial $R^*$ whose roots are
$\eta_i + \frac{b_2}{3b_3}$, i.e.
$$
R^*(x) = 27b_3^2 R\left(x-\frac{b_2}{3b_3}\right) = 27b_3^3 x^3 +
3b_3px +  q \in \ZZ[x],
$$
where the integer coefficients $p,q$ can be explicitly computed.
Since the mutual distances between the roots do not change, compared
with those of $R$, but the leading coefficient of $R^*$ grows by the
factor of $27b_3^2$, we have
\begin{equation}\label{eq5}
\Delta(R^*) = (27b_3^2)^4\Delta(R) \asymp B^8\Delta(R).
\end{equation}
On the other hand, we compute
\begin{equation}\label{eq4}
\Delta(R^*) = 27^2 b_3^6 (-4p^3-27q^2).
\end{equation}

The conclusion is, to compute the lower bound for $\Delta(R)$, we
essentially need to know how small the expression $|4p^3+27q^2|$ can
be. This question was initially investigated by
Hall~\cite{hall_1971} who made a conjecture about its size. It was
later refined a little bit and its current statement is provided in
the Introduction.

%\begin{conjecture}[Hall]
%For any $\epsilon>0$ the inequality $$ 0<|x^3 -
%y^2|<|x|^{1/2-\epsilon}$$ has only finitely many solutions.
%\end{conjecture}

Observe that $|27\cdot 16(4p^3+27q^2)| = |(108 q)^2 - (-12p)^3|$,
Therefore, under the Hall conjecture we can estimate~\eqref{eq4} by
\begin{equation}\label{eq8}
|\Delta(R^*)| \gg_\epsilon B^6 |p|^{1/2-\epsilon},
\end{equation}
where $\epsilon>0$ can be taken arbitrarily small. Also notice that
the roots $\eta^*_1, \eta_2^*, \eta_3^*$ of $R^*$ satisfy
$|\eta_1|\sim \frac23 A, |\eta_i|\sim \frac13 A$, therefore the
coefficient of $R^*$ at $x$ is
\begin{equation}\label{eq10}
|3b_3p| \asymp B^3 A^2\quad \Longleftrightarrow\quad |p|\asymp
(AB)^2.
\end{equation}
In view of~\eqref{eq5} and~\eqref{eq8}, we then get
$$
\Delta(R) \gg_\epsilon \frac{A^{1-2\epsilon}}{B^{1+2\epsilon}}.
$$
The equation~\eqref{sepeq} then implies that $\sep(R)\gg_\epsilon
B^{-2-1/2-\epsilon}A^{-2+1/2-\epsilon}$ or in other words, the pair
$(s,t) = (1/2+\epsilon, 1/2-\epsilon)$ satisfies~\eqref{th1_eq}.
This establishes the claim~\eqref{th2_eq} of Theorem~\ref{th2}. We
conclude that $D_{3,1}$ contains the points $(u,v)$ of the form
$\left(5/2+\epsilon, 5/2+\epsilon\right)$. As $\epsilon$ tends to
zero, this point approaches $(5/2,5/2)$.

Also notice that if the pair $(s_0, t)$ satisfies~\eqref{th1_eq}
then any pair $(s,t)$ satisfies the same inequality as soon as $s\ge
s_0$. Increasing $s$ and keeping $t$ constant, decreases the
parameter $u$ and increases $v$ in~\eqref{eq3}. Then, by denoting $r
= \frac{s}{s+t}$ and changing it between $1/2+\epsilon$ and 1, one
gets that all the pairs
\begin{equation}\label{eq7}
\left(r + \frac{2(1-r)}{1/2-\epsilon}, 2+r\right)
\end{equation}
belong to $D_{3,1}$. This observation provides the elements in the
non-trivial parts of the set $D_{3,1}\cap L_v$, where
$5/2+\epsilon\le v\le 3$. We have thus verified the
claim~\eqref{th2_eq2} of Theorem~\ref{th2}.

If an ineffective version of the Hall conjecture is assumed then we
derive that $|\xi - p/q|\ge c(\epsilon) H^{-u}(\xi)q^{-v}$ for all
$\xi\in\AA_3$, $p/q\in\QQ$ and some $c(\epsilon)>0$. Here, the
parameters $(u,v)$ are given by~\eqref{eq7}. But under this setting,
we can not compute $c(\epsilon)$, we only know that it exists. One
can improve the situation a little bit by effectively constructing
the constant $c>0$ and proving the inequality $|\xi-p/q|\ge
cH^{-u}(\xi)q^{-v}$ for all but finitely many pairs $(\xi,p/q)\in
\AA_3\times \QQ$. Indeed, we can replace the sign $\gg_\epsilon$ by
$\ge$ in~\eqref{eq8} for all but finitely many values $p$. Since
$p\asymp AB$, there are only finitely many pairs $A,B$ failing that
inequality. But in view of $A\asymp q^{v-2}H^u$, this implies the
finiteness of the set of pairs ($\xi,p/q$) that may produce $A,B$
that fail an updated inequality~\eqref{eq8}. For the remaining pairs
$(\xi,p/q)$ all the inequalities $\gg_\epsilon$ in the proof can be
replaced by $\gg$ and all the implied absolute constants can be
explicitly computed. This establishes the claim~\eqref{th2_eq3} of
Theorem~\ref{th2}.

Finally, we note that if one assumes the effective version of the
Hall conjecture then the inequality~\eqref{main_eq} can be made
effective as well. More precisely, assume that for some $\epsilon>0$
one can provide a constant $c>0$ such that for all $x,y\in\ZZ$ if
$x^3\neq y^2$ then
$$
|x^3 - y^2|>cx^{\frac12-\epsilon}.
$$
Then one can explicitly compute another constant $c^*$ such that
$$
|\xi - p/q| > \frac{c^*}{H^u(\xi)q^v}
$$
with $(u,v)$ given by~\eqref{eq7}. Theorem~\ref{th2} is now fully
proven.

Now we show that assumption of the more general abc-conjecture leads
to better estimates on $(u,v)$. Recall that from~\eqref{eq5}
 and~\eqref{eq4} we have that $-4p^3 - 27q^2$ is a multiple of
$b_3^2$. We write it as $b_3^2d$ where $d\in\ZZ$. Now we denote
$g:=\gcd(-4p^3,27q^2)$,
$$
a := \frac{b_3^2 d}{g},\quad b :=\frac{27q^2}{g},\quad c :=
\frac{-4p^3}{g}
$$
and get that $a,b,c$ are coprime with $a+b=c$. Notice that
$$
\rad(a) \le \frac{b_3d}{\sqrt{g}},\quad \rad(b)\le
\frac{q}{\sqrt{g}}\quad\mbox{and}\quad \rad(c)\le
\frac{p}{\sqrt[3]{g}}.
$$

Assume that $a\le p$ then $p^3\asymp q^2$ and in view
of~\eqref{eq10},
$$
p\asymp A^2B^2,\quad q\asymp A^3B^3,\quad \max\{|a|,|b|,|c|\} \asymp
A^6B^6g^{-1}, \rad(abc) \le b_3pqd\asymp A^5B^6dg^{-4/3}.
$$
Assuming the abc-conjecture, we get that the inequality
$$
\max\{|a|,|b|,|c|\}\ll_\delta \rad(abc)^{1+\delta}\quad
\Longrightarrow\quad A^6B^6\ll_\delta (A^5B^6d)^{1+\delta}
$$
is satisfied for all triples $b_3,p,q$. by choosing $\delta =
\frac{\epsilon/3}{1-\epsilon/3}$, where $\epsilon$ is an arbitrary
positive number and can be as small as we wish, we derive
$$
A^5B^6d \gg_\epsilon A^6B^6 (AB)^{-2\epsilon}\quad
\Longleftrightarrow\quad d\gg_\epsilon
\frac{A^{1-2\epsilon}}{B^{2\epsilon}}\quad \Longleftrightarrow\quad
\Delta(R^*)\gg_\epsilon B^{8-2\epsilon}A^{1-2\epsilon}.
$$
Finally, in view of~\eqref{eq5}, we derive
$$
\Delta(R)\gg_\epsilon \frac{A^{1-2\epsilon}}{B^{2\epsilon}}.
$$

For the remaining part of the proof, we proceed analogously to the
case under the assumption of the Hall conjecture.
Equation~\eqref{sepeq} provides $\sep(R)\gg_\epsilon
B^{-2-\epsilon}A^{-2+1/2 - \epsilon}$ or in other words, the pair
$(s,t) = (\epsilon, 1/2-\epsilon)$ satisfies~\eqref{th1_eq}.
Therefore $D_{3,1}$ contains the points of the form $(4+2\epsilon,
2+2\epsilon)$ which approaches $(4,2)$ as $\epsilon\to 0$. Finally,
by considering pairs $(s,t) = (s,1/2-\epsilon))$ with $s\ge
\epsilon$ we derive that $D_{3,1}$ contains the pairs of the
form~\eqref{eq7} for all $r\in[2\epsilon,1]$.

\section{Application to Thue equations}\label{sec6}

The relation between rational approximations to irrational algebraic
numbers and Thue equations is well know. We provide it here for
completeness. Write
$$
|F_\va(p,q)| = q^3|P_\va(p/q)|,
$$
where $P_\va\in\ZZ[x]$ is the cubic polynomial whose coefficient
vector is $\va$. Then, by letting $\xi$ to be the closest root of
$P_\va$ to $p/q$, we derive
$$
|F_\va(p,q)| \asymp  q^3\left|\xi - \frac{p}{q}\right|\cdot
|a_3(\xi-\xi_2)(\xi-\xi_3)|.
$$
Here we follow the notation from Section~\ref{sec2}, i.e.
$\xi_2,\xi_3$ are the remaining roots of $P_\va$. This time, we need
to provide the lower bound for the right hand side.

If, among the three roots of $P_\va$, $\xi_2$ and $\xi_3$ are
closest to each other, then we have
$$
|a_3|(\xi-\xi_2)(\xi-\xi_3)|\ge |a_3
\sqrt{|(\xi-\xi_2)(\xi-\xi_3)(\xi_2-\xi_3)|}\cdot \sqrt{|\xi-\xi_3|}
$$
Otherwise, without loss of generality we can assume that $\xi$ and
$\xi_2$ are the closest roots of $P_\va$. In this case we have
$|\xi-\xi_3| \asymp \sqrt{|(\xi_2-\xi_3)(\xi-\xi_3)|}$. In both
cases, we derive
$$
|a_3(\xi-\xi_2)(\xi-\xi_3)| \gg \Delta^{1/4}(P_\va)\sep^{1/2}(P_\va)
\gg \frac{1}{H(P_\va)} = \frac{1}{||\va||_\infty}.
$$
The last inequality follows from $\sep(P_\va)\gg H^{-2}(P_\va)$, due
to Mahler~\cite{mahler_1964}.

The conclusion is that one always has
$$
|F_\va(p,q)| \gg \frac{q^3}{||\va||_\infty}\cdot
\left|\xi-\frac{p}{q}\right|.
$$
If $\xi$ is quadratic irrational or a rational number distinct from
$p/q$ then Liouville's inequality implies $|\xi - p/q|\gg
q^{-2}H^{-1}(\xi)$. Finally, we apply Gelfond's lemma which states
that for all polynomials $P, Q\in \RR[x]$, $H(P)H(Q) \asymp H(PQ)$.
For the proofs of Liouville's inequality and Gelfond's lemma we
refer the reader to~\cite[Appendix A]{bugeaud_2004}. They imply that
$H(\xi)\ll ||\va||_\infty$ and then the claim of
Corollary~\ref{corl} follows immediately. On the other hand, if
$\xi$ is cubic irrational then Corollary~\ref{corl} follows from the
claim~\eqref{th2_eq3} of Theorem~\ref{th2} with $v = 5/2 + \epsilon$
under the assumption of the Hall conjecture; or with $v =
2+2\epsilon$ under the assumption of the abc-conjecture.

\section{Family of cubic irrationals}\label{sec4}

Here we prove Theorem~\ref{th3} by constructing an infinite family
of cubic irrationals $\xi_n$ and corresponding rational
approximations $p_n/q_n$ such that
$$
\left|\xi_n - \frac{p_n}{q_n}\right| \ll \frac{1}{
H(\xi_n)^{4-3v}q_n^{2+v}}
$$
for all $1\le v\le 2$.

Let $v_n/u_n$ be the $n$'th convergent to $\sqrt{2}$. Since
$\sqrt{2} = [1;\overline{2}]$, easy computations give
$$
v_0=1;\quad v_1=3;\quad v_{n+1} = 2v_n + v_{n-1},\; \forall n\ge 2;
$$
$$
u_0=1;\quad u_1=2;\quad u_{n+1} = 2u_n + u_{n-1},\; \forall n\ge 2.
$$
From basic properties of convergents together with the theory of
Pell's equations we get that
\begin{equation}\label{eq6}
2u_n^2 - v_n^2 = (-1)^n;\quad u_{n-1}v_n - v_{n-1}u_n = (-1)^n.
\end{equation}
Consider the following family of cubic polynomials:
\begin{equation}\label{eq15}
P_n(x):= u_nx^3 + (4u_n+v_n)x^2 - 2(4v_n - v_n)x -2v_n.
\end{equation}
Enumerate the roots of $P_n(x)$ by $\xi_{1,n}, \xi_{2,n},
\xi_{3,n}$. One can check that the ratios between the absolute
values of the coefficients of these polynomials converge, more
precisely,
$$
\lim_{n\to\infty} \frac{5u_n+u_{n-1}}{u_n} =
5+\frac{1}{\sqrt{2}+1};\; \lim_{n\to\infty}\frac{6v_n-2u_{n-1}}{u_n}
= 6\sqrt{2} - \frac{2}{\sqrt{2}+1};\lim_{n\to\infty}
\frac{2v_n}{u_n} = 2\sqrt{2}.
$$
Indeed, these limits follow from the one of the properties of the
convergents, namely that $$\frac{v_n}{u_n},\;
\frac{v_{n-1}}{u_{n-1}}\to \sqrt{2},\quad\frac{u_{n}}{u_{n-1}} =
[\underbrace{2;2,2,\ldots,2}_{n \mbox{ times}}]\; \to\; 1+\sqrt{2}.
$$
Therefore the sequences $\xi_{i,n}$, $1\le i\le 3$ converge to the
roots of the polynomial
$$
P_\infty(x) = x^3 + (4+\sqrt{2})x^2 - (4\sqrt{2}+2)x - 2\sqrt{2}.
$$
An easy exercise is to check that one of them is $\sqrt{2}$. Without
loss of generality, one may assume that $\xi_{1,n}\to \sqrt{2}$ as
$n$ tends to infinity. Notice that these roots stay away from each
other, but are not distanced too far apart. In other words, we have
that $|\xi_{i,n}-\xi_{j,n}|\asymp 1$ for all $n\in\NN$ and $1\le
i\neq j\le 3$.

Consider the following rational number:
\begin{equation}\label{eq12}
\frac{p_n}{q_n}:=
\frac{4u_n^3+16u_n^2v_n+14u_nv_n^2+4v_n^3}{8u_n^3+14u_n^2v_n+8u_nv_n^2+v_n^3}.
\end{equation}
Notice that $H(P_n)\asymp v_n$ while $q_n\asymp v_n^3$. Also, one
can easily check that
$$
\lim_{n\to\infty} \frac{p_n}{q_n} = \frac{(4+16\sqrt{2} + 14\cdot 2
+ 4\cdot2\sqrt{2})u_n^3}{(8+14\sqrt{2}+8\cdot 2+ 2\sqrt{2})u_n^3} =
\frac{4+3\sqrt{2}}{3+2\sqrt{2}} = \sqrt{2}.
$$
Therefore, for $n$ bigger than some value $n_0$, $\xi_{1,n}$ is the
closest root to $p_n/q_n$. A brave reader can use
equations~\eqref{eq6} to verify that
\begin{equation}\label{eq14}
q_n^3 P_n(p_n/q_n) = 2.
\end{equation}
Alternatively, one can argue as follows. In view of~\eqref{eq15}
and~\eqref{eq12}, $q_n^3P_n(p_n/q_n)$ is a homogeneous polynomial of
degree 10 in $u_n$ and $v_n$. Next, in view of the recurrent
formulae for $u_n$ and $v_n$, these quantities are elements of the
linearly recurrent sequences with characteristic polynomial
$x^2-2x-1$ and therefore they can be expressed as
$$
u_n = \lambda_1 \alpha_1^n + \lambda_2\alpha_2^n, \quad v_n =
\tau_1\alpha_1^n + \tau_2\alpha_2^n,
$$
where $\alpha_1,\alpha_2$ are the roots of the characteristic
polynomial and $\lambda_1,\lambda_2,\tau_1,\tau_2$ are constants
that can be explicitly computed. We derive that $q_n^3P_n(p_n/q_n)$
can be written as a linear combination of the quantities
\begin{equation}\label{eq13}
(\alpha_1^{10})^n,\; (\alpha_1^9\alpha_2)^n,\; \ldots,\;
(\alpha_2^{10})^n,
\end{equation}
where all the coefficients can be explicitly computed. This implies
that $q_n^3P_n(p_n/q_n)$ is a linearly recurrent sequence of degree
11, whose characteristic polynomial has roots from~\eqref{eq13}.
Finally, with help of a computer package such as Mathematica, one
can verify~\eqref{eq14} for all $n$ between 1 and 11, hence our
sequence must be constant for all $n\in\NN$.

Equation~\eqref{eq14} can be rewritten as
$$
2 = q_n^3 u_n \left(\xi_{1,n} -
\frac{p_n}{q_n}\right)\left(\xi_{2,n} -
\frac{p_n}{q_n}\right)\left(\xi_{3,n} - \frac{p_n}{q_n}\right)
\asymp q_n^3v_n\left|\xi_{1,n} -\frac{p_n}{q_n}\right|.
$$
This immediately leads to
$$
\left|\xi_{1,n} - \frac{p_n}{q_n}\right| \ll \frac{1}{q_n^3v_n}
\asymp \frac{1}{q_n^{2+v}H^{4-3v}(P_n)}.
$$

In fact, an enthusiastic reader may want to verify that the
continued fraction of the root $\xi_{1,n}$ starts with
$$
\xi_{1,n} = [1;\underbrace{2,2,\ldots,2}_{n}4,
\underbrace{2,2,\ldots,2}_{2n+1}, 3,  A_n, 1,1,
\underbrace{2,\ldots,2}_{2n+1},
1,1,1,\underbrace{2,\ldots,2}_{2n+1},1,1,\ldots]
$$
where
$$
A_n = 2 + 4u_nv_n (14u_n^2 + 20 u_nv_n + 7v_n^2).
$$
We do not verify this claim here, but just want to mention how
unusual the continued fractions of some cubic irrationals may look
like. Numerical computations suggest that the partial quotients of
$\xi_{1,n}$ with higher indices seem to look pretty random.
Therefore, despite the very unusual pattern at the beginning, the
distribution of partial quotients of $\xi_{1,n}$ may still be the
same as for a generic real number.

\section{The case of function fields}\label{sec5}

Let $\FF$ be a field. Consider the function field $\KK:=
\FF((t^{-1}))$ of formal Laurent series in $t^{-1}$ and define the
norm in this space as follows. Let $$x = \sum_{k=-d}^\infty
a_kt^{-k} \in \FF((t^{-1}))$$ be such that $a_{-d}\neq 0$. Then
$|x|:= 2^d$. A more standard approach is to define it as $|x|:=p^d$
where $p$ is the characteristic of $\FF$. However, it does not work
for fields of characteristic zero
 and for our purposes the exact base in the absolute value is
 not important.

 One can easily see that $\FF[t]$ and $\FF(t)$ are both
 subsets of $\FF((t^{-1}))$. Equipped with the norm, we can
 construct a theory of Diophantine approximation similar to that in
 $\RR$, where $\FF[t]$ and $\FF(t)$ play the role of integer and rational
 numbers respectively.

For $d\in \NN$ we define $\AA_d(\KK)$ to be the set of algebraic
series $x\in\KK$ whose degree is exactly $d$, i.e. the series that
satisfy the equation $P(x) = 0$ where $P\in \FF[t][x]$ has degree
$d$ and is irreducible as a polynomial in $x$. One can make this
polynomial minimal by requiring its coefficients to be coprime
polynomials in $\FF[t]$. By the height of $\alpha\in \AA_d(\KK)$ we
define the maximum of the norms of the coefficients of the minimal
polynomial of $\alpha$. Using this notion, one can define the sets
$D_{d,r}(\KK)$ analogously to the sets $D_{d,r}$: we say that
$(u,v)\in D_{d,r}(\KK)$ if there exists an absolute constant
$c_{d,r}$ such that for all $\xi\in \AA_d(\KK)$ and $\alpha\in
\AA_r(\KK)$ one has
$$
|\xi-\alpha| > \frac{c_{d,r}}{H(\xi)^{u}H(\alpha)^{v}}.
$$

It is known that an analogue of Liouville's inequality also holds in
this setting. Namely, for algebraic elements $\xi, \alpha\in\KK$ of
degrees $d$ and $r$, one has $|\xi-\alpha| \gg
H(\xi)^{-r}H(\alpha)^{-d}$. Consequently,
$$
B_{r,d}:= \{(u,v)\in\RR: u\ge r, v\ge d\} \subset D_{d,r}(\KK).
$$
For the case $r=1$ (i.e, when $\alpha$ is rational), this inequality
was established by Mahler~\cite{mahler_1949}. For higher values of
$r$, one can verify the inequality by following the same arguments
as in~\cite[Theorem A.1]{bugeaud_2004} adapted to algebraic elements
in function fields.

It appears that the remaining part $D_{d,r}(\KK)\setminus B_{r,d}$
is much better understood than in the real case. Its shape depends
in an essential way on whether the characteristic of $\FF$ is zero
or positive. From now on, we restrict our attention to the case
$r=1$.

Suppose that $\FF$ has characteristic $p>0$.
Mahler~\cite{mahler_1949} first observed that the root $\alpha$ of
the equation $t\alpha^d - t\alpha - 1=0$, where $d$ is a power of
$p$, admits infinitely many rational approximations $P/Q\in\FF(t)$
such that
$$
|\alpha - P/Q|\ll |Q|^{-d}.
$$
This shows that any $(u,v)\in D_{d,1}$ must satisfy $v\ge d$. Later,
Osgood~\cite{osgood_1975} constructed examples of algebraic series
of arbitrary degree for which Liouville's inequality is sharp. This
yields a full description of $D_{d,1}(\KK)$:
$$
D_{d,1}(\KK) = B_{1,d}.
$$

If $\FF$ has characteristic zero the situation is completely
different. As shown by Osgood~\cite{osgood_1973}, cubic irrational
series in this setting are badly approximable. Moreover, one can
prove that $D_{3,1}(\KK)\cap L_2$ is non-empty. For $d>3$ his result
is weaker but still yields $D_{d,1}(\KK)\cap L_\lambda\neq\emptyset$
for $\lambda \gg (\log d)^2$. From now on, we focus on the set
$D_{3,1}(\KK)$ and use the arguments of that paper to obtain an
estimate for $u$ such that $(u,2)\in D_{3,1}(\KK)$.

By replacing $\alpha$ with $\alpha^{-1}$ if needed, we can assume
that $|\alpha|\le 1$. Next, by replacing $\alpha$ with $\alpha+1$ if
needed, we can make sure that $|\alpha|=1$. Next, as was shown
in~\cite{osgood_1973}, every cubic irrational $\alpha\in \KK$ is
also a solution of a Riccatti equation
$$
Dx' = A+Bx+Cx^2,\quad A,B,C,D\in \FF[t].
$$
Now, let $\alpha = \frac{P}{Q} + \delta$ for some $\frac{P}{Q}\in
\FF(t)$. Substituting this into the Riccatti equation gives
$$
D\left(\frac{P}{Q}\right)' - A - B\frac{P}{Q} - C\frac{P^2}{Q^2} =
B\delta + 2C\frac{P}{Q}\delta + C\delta^2 - D\delta'.
$$
Note that the left hand side is a rational function with the
denominator $Q^2$. It is also shown in~\cite{osgood_1973} to be
nonzero. Therefore its norm is at least $|Q|^{-2}$. Estimating the
norm of the right hand side gives
$$
|Q|^{-2} \le \max\left\{|B||\delta|, |C||\delta|,
|D|\frac{|\delta|}{2}\right\}.
$$
We conclude that $|\delta|\ge (\max\{|B|, |C|, |D|/2\} |Q|^2)^{-1}$.

Finally, the precise expression~\cite[Proposition 1]{badziahin_2024}
for the coefficients $A,B,C,D$ in terms of the coefficients of the
minimal polynomial $P_\alpha$ gives $\max\{|A|, |B|, |C|, |D|\}\le
|H(\alpha)|^4$. We conclude
$$
\left|\alpha - \frac{p}{q}\right|\ge \frac{1}{H(\alpha)^4 H(P/Q)^2},
$$
i.e. $(4,2)\in D_{3,1}(\KK)$.

On the other hand, there are examples of cubic series $\alpha$ whose
continued fractions have infinitely many partial quotients $a_n$
such that $|a_n|\ge 3H(\alpha)$. For example~\cite{badziahin_2024},
the solution of $3x^3 - 3tx^2 - 3cx + ct=0$ where $c\in \QQ$ is a
parameter, has a continued fraction
\begin{equation}\label{no5}
\bigk\left[ \begin{array}{l}  \\ t \end{array} \overline{\begin{array}{llll}2(3k+1)c&(6k+1)c&2(3k+2)c^2& (6k-1)c^2\\
3(4k+1)t&t&3(4k+3)t(t^2+2c)&t
\end{array}\hspace{-1ex}}\;\; \right].
\end{equation}
The partial quotients here are overlined to indicate that the rule
for them is periodic, while the partial quotients themselves are
not. The term $k$ in the notation indicates the periodic template's
number where we start counting from $k=0$. One can replace $t$ in
the minimal polynomial and the continued fraction by any
$R(t)\in\QQ[t]$ of arbitrary degree. Then, by considering $P/Q$ to
be $(4i+2)$'th convergents of $\alpha$, we find infinitely many
pairs $(\alpha, P/Q)$ such that
$$
\left|\alpha - \frac{P}{Q}\right|\le \frac{1}{H(\alpha)^3 H(P/Q)^2}
$$
where $H(\alpha)$ and $H(P/Q)$ can both be made arbitrarily large.
Hence any $(u,2)\in D_{3,1}(\KK)$ must satisfy $u\ge 3$.

An interesting problem is to find the minimal value of $u$ for
which$(u,2)\in D_{3,1}(\KK)$. From the above discussion, this value
must satisfy $3\le u\le 4$. It would also be great to describe, or
at least estimate, the shape of $D_{3,1}(\KK)$ for $2<v<3$. However,
since this paper focuses instead on the different set
$D_{3,1}\subset \RR^2$, we leave these questions for future
investigation. But the main takeaway is that the set $D_{3,1}(\KK)$
is nontrivial for all $2\le v\le 3$ when $\FF$ has characteristic
zero. This stands in sharp contrast to the case of positive
characteristic, where any $(u,v)\in D_{3,1}(\KK)$ must satisfy $v\ge
3$.

In Diophantine approximation it is quite common for results over
$\RR$ and over function fields to exhibit strong similarities, if
not direct analogies. However, the above arguments show that this
parallel can not be expected when describing the shape of $D_{d,r}$,
since the sets $D_{d,r}(\KK)$ behave very different depending on the
characteristic of $\FF$. Nevertheless, Theorem~\ref{th5} suggests
that
%if we consider $\KK
%= \QQ((t^{-1}))$ and specialise a given algebraic $\alpha\in\KK$ and
%its best approximant $p(t)/q(t)$ at some integer value $t$, we
%derive an algebraic $\alpha(t)\in\RR$ and a rational
%$p(t)/q(t)\in\QQ$. In many cases, this rational number appears to
%approximate $\alpha(t)$ well. Therefore, the problems about
%$D_{d,r}$ in $\RR^2$ and $D_{d,r}(\KK)$ for $\KK=\QQ((t^{-1}))$ are
%somehow related. This heuristic and a bit vague argument suggests
$D_{d,r}$ may resemble $D_{d,r}(\KK)$ more closely in the zero
characteristic setting.

\bigskip
\noindent Dzmitry Badziahin\\ \noindent The University of Sydney\\
\noindent Camperdown 2006, NSW (Australia)\\
\noindent {\tt dzmitry.badziahin@sydney.edu.au}

\section*{Statements and Declarations}

The author declares that no funds, grants, or other support were
received during the preparation of this manuscript.

\noindent The authors have no relevant financial or non-financial
interests to disclose.

\end{document}